\documentclass{amsart}
\usepackage{epsfig}
\newenvironment{theorem}{\noindent\normalsize {\sc Theorem.}}{}
\newenvironment{conjecture}{\noindent\normalsize {\sc Conjecture.}}{}
\newenvironment{lemma}{\noindent\normalsize {\sc Lemma.}}{}
\newenvironment{prop}{\noindent\normalsize {\sc Proposition.}}{}

\def\E{{\mathcal E}}

\def\Hy{{\mathcal H}}
\def\cO{{\mathcal O}}

\def\Z{{\mathbb Z}}
\def\Q{{\mathbb Q}}
\def\R{{\mathbb R}}

\def\A{{\mathbb A}}
\def\scS{{\it\mathcal S}}
\def\dim{{\rm dim\:}}
\def\vol{{\rm vol\ }}

\begin{document}
\title{On Volumes of Arithmetic Quotients of $SO(1,n)$.}
\author{Mikhail Belolipetsky}
\address{Sobolev Institute of Mathematics, Koptyuga 4, 630090 Novosibirsk, Russia
\newline\phantom{iii}
Max Planck Institute of Mathematics, Vivatsgasse 7, 53111 Bonn, Germany}

\email{mbel@math.nsc.ru}
\date{}

\begin{abstract}
We apply G.~Prasad's volume formula for the arithmetic quotients
of semi-simple groups and Bruhat-Tits theory to study the
covolumes of arithmetic subgroups of $SO(1,n)$. As a result we
prove that for any even dimension $n$ there exists a unique
compact arithmetic hyperbolic $n$-orbifold of the smallest
volume. We give a formula for the Euler-Poincar\'e characteristic
of the orbifolds and present an explicit description of their fundamental
groups as the stabilizers of certain lattices in quadratic spaces.
We also study hyperbolic $4$-manifolds defined arithmetically and
obtain a number theoretical characterization of the smallest compact
arithmetic $4$-manifold.
\end{abstract}

\maketitle

\section{Introduction}

In this article we consider the problem of determining the smallest hyperbolic manifolds
and orbifolds defined arithmetically. This problem has a long history which goes back to
Klein and Hurwitz. Its solution for the hyperbolic dimension~$2$ was known to
Hurwitz which allowed him to write down his famous bound for the order of the automorphisms
group of a Riemann surface. The first extremal example for the bound is the Klein quartic.
Many interesting facts about this classical subject and far reaching generalizations can be
found in the book \cite{Levi}. For the dimension $3$ the problem is also completely solved
but the results are quite recent \cite{CF1}, \cite{CF2}. For the higher dimensions very
little is known.

Probably the most interesting case among the dimensions higher than $3$ is in dimension
$4$. Recently it has attracted particular attention due to a possible application in cosmology:
closed orientable hyperbolic $4$-manifolds arise as the doubles of the real tunnelling
geometries if the cosmological constant is assumed to be negative. In this context
there are physical arguments in favor of using the
smallest volume orientable hyperbolic $4$-manifold as a
model of the Lorentzian spacetime \cite{Gib}. In view of the known facts and
conjectures for the small dimensions it is quite natural to look for the smallest manifold
among the arithmetic ones.

In this article we apply G.~Prasad's volume formula for the arithmetic quotients of semi-simple
groups and Bruhat-Tits theory to investigate the particular case of the group $SO(1,n)$ whose
symmetric space is the hyperbolic $n$-space.
Since our primary interest lies in dimension $4$, at some point we restrict our attention
to even dimensions and in the final section even more restrictively, we consider only the
$SO(1,4)$-case. The main results are given in Theorems~\ref{thm1} and \ref{thm2}.

The first theorem says that the smallest compact arithmetic orbifold in any even dimension
greater or equal than $4$ is unique and defined over the field $\Q[\sqrt{5}]$, it also provides
an explicit description for the orbifold and a formula for its volume. Let us remark that the
quadratic number field $\Q[\sqrt{5}]$ has the interesting property that its fundamental unit
$\epsilon = (1 + \sqrt{5})/2$ is the ``golden section'' unit which was already known to
Greek mathematicians.
Since for the dimension $2$ the situation is different (the field of definition
of the smallest hyperbolic $2$-orbifold is $\Q[\cos(2\pi/7)]$) this result was a little
unexpected for us.

The problem of determining the smallest arithmetic manifold is much more delicate than that
for the orbifolds. Here we are currently able to present only partial results and only in
dimension $4$. Still our Theorem~\ref{thm2} gives an explicit classification of
all the possible candidates and essentially reduces the problem of finding the smallest
arithmetic $4$-manifold to an extensive computation, which we hope is practically
possible.
\medskip


I would like to thank Professor G.~Harder for helpful discussions and encouragement.
I would also like to thank Wee Teck Gan for the helpful email correspondence.
Finally, I would like to thank Professor S.~Lang for reading an early version of this
paper and for his suggestions about the presentation.
While working on this paper I was visiting Max-Planck-Institut f\"ur Mathematik in Bonn,
I appreciate the hospitality and financial support from MPIM.

\section{Arithmetic subgroups}

\subsection{}
Hyperbolic $n$-space can be obtained as a symmetric space associated to the orthogonal group $G$
of type $(1,n)$:
$$ \Hy^n = G/K_G = SO(1,n)^o/SO(n)$$
($K_G$ denotes a maximal compact subgroup of $G$).
This way the connected component of identity of $SO(1,n)$ acts as a group of
isometries of $\Hy^n$. A discrete subgroup $\Gamma$ of $G$ defines a locally
symmetric space $X = \Gamma\backslash G/K_G$ which in our case will be an orientable hyperbolic
orbifold or manifold if $\Gamma$ is torsion-free. We will be interested in hyperbolic orbifolds
which arise from arithmetic subgroups of $G$.

Let $G$ be a connected semi-simple Lie group, $H/k$ is a semi-simple algebraic group defined
over a number field $k$ and $\phi: H(k\otimes\R) \to G$ is a surjective homomorphism with a
compact kernel. We consider $H$ as a $k$-subgroup of $GL(n)$ for $n$ big enough and define a
subgroup $\Lambda$ of $H(k)$ to be {\it arithmetic} if it is commensurable with the subgroup of
$k$-integral points $H(k)\cap GL(n,\cO_k)$, that is, the intersection $\Lambda\cap GL(n,\cO_k)$
is of finite index in both $\Lambda$ and $H(k)\cap GL(n,\cO_k)$.
The subgroups of $G$ which are commensurable with $\phi(\Lambda)$ are called {\it arithmetic
subgroups} of $G$ defined over the field $k$.
It can be shown that the notion of arithmeticity does not depend on a particular choice of
the $k$-embedding of $H$ into $GL(n)$.

We call an orbifold or manifold $X = \Gamma\backslash G/K_G$ {\it arithmetic} if $\Gamma$
is an arithmetic subgroup of $G$, and we say that $X$ is {\it defined over $k$} if $k$ is the
field of definition of~$\Gamma$.

Arithmetic subgroups of the orthogonal groups can be constructed as follows.
Let now $k$ be a totally real algebraic number field with the ring of integers $\cO$ and
let $f$ be a quadratic form of type $(1,n)$ with the coefficients in $k$ such that
for any non-identity embedding $\sigma:k\to\R$ the conjugate form $f^\sigma$ is
positive definite (such an $f$ is called {\it admissible form}). Then given an
$\cO$-integral lattice $L$ in $k^{n+1}$ the group $\Gamma = G_L =
\{\gamma\in G \cong SO(f)^o \mid L\gamma = L\}$ is an arithmetic subgroup
of $G$ defined over $k$.

It can be shown that for even $n$ this construction gives all arithmetic subgroups
of $G  = SO(1,n)^o$ up to commensurability and conjugation in $G$. For odd $n$ there
is also another construction related to quaternion algebras, and for $n = 7$ there is a
special type of arithmetic subgroups related to the Cayley algebra.

\subsection{}\label{max1}
Looking for the hyperbolic orbifolds and manifolds of the smallest volume we will be interested
in the {\it maximal arithmetic subgroups} of $G$. The maximal arithmetic subgroups can be
effectively classified in terms of the arithmetic data and the local structure of $G$.
In order to discuss the classification picture we will give some more definitions.

Let $G/k$ be a connected semi-simple algebraic group defined over a number field $k$, and
let $V_f$ (resp. $V_\infty$) denote the set of finite (resp. infinite) places of $k$.
By \cite{BT} for a local place $v\in V_f$ the group $G(k_v)$ is endowed with the
structure of Tits system of affine type $(G(k_v), B_v, N_v, \Delta_v)$.
A subgroup $I_v\subset G(k_v)$ is called {\it Iwahori
subgroup} if it is conjugate to $B_v$. A subgroup $P_v\subset G(k_v)$ which contains an
Iwahori subgroup is called {\it parahoric}.
A collection $P = (P_v)_{v\in V_f}$ of parahoric subgroups $P_v$ is said to
be {\it coherent} if $\prod_{v\in V_\infty}G(k_v)\cdot\prod_{v\in V_f} P_v$ is an open subgroup
of the ad\`ele group $G(\A)$. A coherent collection of parahoric subgroups
$P = (P_v)_{v\in V_f}$ defines an arithmetic subgroup $\Lambda = G(k)\bigcap\prod_{v\in V_f}
P_v$ of $G(k)$ which will be called the {\it principal arithmetic subgroup} determined by $P$.
We will also call the corresponding subgroups of the Lie group $G$ principal arithmetic
subgroups.

If the group $G/k$ is simply connected and adjoint then the maximal arithmetic subgroups
of $G$ are exactly the principal arithmetic subgroups defined by coherent collections
of maximal parahoric subgroups.  For the other forms the
situation becomes more complicated, but still it is true that any maximal arithmetic subgroup
is a normalizer in $G$ of some principal arithmetic subgroup \cite{Plat}. The problem of
classification of the principal arithmetic subgroups which give rise to the maximal
subgroups was studied in \cite{CR} where, in particular, a criterion for the groups
of type $B_r$ (this is the type of $SO(1,2r)$) is given. However, the criterion
of Ryzhkov and Chernousov is a little subtle: it provides explicit conditions
on the collections of parahoric subgroups but it does not always
guarantee the existence of the global subgroup with the prescribed local properties.
Let us consider this more carefully.

\subsection{}\label{par_cond}
Let $\Lambda$ be a principal arithmetic subgroup of $G/k$ defined by
$\prod_{v\in V_f}P_v\subset G(\A_f) = \prod'_{v\in V_f} G(k_v)$ (where $\Pi'$
denotes the restricted product with respect to $G(\cO_v)$). For each place $v$ the
type of $P_v$ depends on the splitting type of $G(k_v)$. We claim that there is a
natural restriction on the possible splitting types of $G(k_v)$ which, in turn,
implies a restriction on types of $P_v$.

In \cite{K} to any reductive group $G/k$ Kottwitz assigned an invariant $\epsilon(G)\in\{\pm 1\}$,
which can be computed explicitly for $G$ over the completions of $k$ and for which
the product formula holds. Thus, for $G$ over a nonarchimedean local field
$$\epsilon(G(k_v)) = (-1)^{r(G_{qs}) - r(G)}$$
($G_{qs}$ denotes the quasi-split inner
form of $G$ and $r(G)$ is the $k_v$-rank of the derived group of $G$), and for
the archimedean places
$$\epsilon(G(k_v)) = (-1)^{q(G_{qs}) - q(G)}$$
($q(G)$ is a half of the dimension of the symmetric space attached to $G(k_v)$). From
the product formula for $\epsilon(G)$ we immediately obtain that the total number of
places for which $\epsilon = -1$ is even. This is what we call the {\it parity
condition} on the number of nonsplit places. Let us see what does it mean for
our semi-simple groups of type~$B_r$.

There are two forms $B_r$ and $^2B_r$ of type $B_r$ over a nonarchimedean local field, the first form
is split and the second is a non quasi-split form \cite{Tits}. In the first case the $k_v$-rank
$r(G) = r$ and in the second case $r(G) = r-1$, and also always $r(G_{qs}) = r$.
So for $v\in V_f$ we have:
$$
\epsilon(G(k_v)) = \left\{
\begin{array}{r@{\quad}l}
1  & {\rm if\ } G {\rm \ is\ split\ over\ } k_v,\\
-1 & {\rm otherwise.}\\
\end{array} \right.
$$

Over the archimedean places of $k$ by the admissibility condition $G(k_v)$ is compact for all $v\in V_\infty$
except at one place, say $v = Id$, which implies:
$$q(G(k_v)) = 0\ {\rm for}\ v\neq Id,$$
$$q(G(k_{Id})) = r.$$
For the quasi-split form we have:
\begin{eqnarray*}
q(G_{qs}) & = & \frac12 (\dim (G_{qs}) - \dim (K_{G_{qs}})) \\
& = & \frac12 (\dim (SO(r+1,r)) - \dim (S(O(r+1)\times O(r)))) = (r^2 + r)/2.
\end{eqnarray*}
So for the place over which $G$ is non-compact ($v = Id$):
$$
\epsilon(G(k_v)) = \left\{
\begin{array}{r@{\quad {\rm if}\quad}l}
1  & r\equiv 0,1\ ({\rm mod}\ 4),\\
-1 & r\equiv 2,3\ ({\rm mod}\ 4);\\
\end{array} \right.
$$
and for all the other infinite places:
$$
\epsilon(G(k_v)) = \left\{
\begin{array}{r@{\quad {\rm if}\quad}l}
1  & r\equiv 0,3\ ({\rm mod}\ 4),\\
-1 & r\equiv 1,2\ ({\rm mod}\ 4).\\
\end{array} \right.
$$
This implies, for example, that over a totally real quadratic field
the number of nonarchimedean places over which $G$ does not split is odd for odd $r$
and even for $r$ even. So we have a parity condition on the number of places over
which $G$ has type $^2B_r$. As we already remarked, this gives a restriction on the
possible types of the collections of parahoric subgroups of $G$. This kind of restriction
can not be seen in \cite{CR}, but it appears to be important for the applications.

It can be checked using the stabilizers of lattices at least for the orthogonal
groups that this condition is also sufficient for the existence of the prescribed
collections of parahoric subgroups.

\subsection{}
From the general theory of arithmetic subgroups of semi-simple Lie groups it follows
that any arithmetic subgroup $\Gamma$ is a discrete subgroup of $G$ and the volume of
$\Gamma\backslash G/K_{G}$ is finite. It is also known that
$\Gamma\subset SO(1,2r)^o$ is cocompact if and only if the corresponding
quadratic form $f$ does not represent zero non-trivially over $k$, which for $n = 2r\geq 4$
means that $\Gamma$ is non-cocompact if and only if it is defined over $\Q$.
We are going to investigate the volumes of $\Gamma\backslash\Hy^n$.

\medskip

For the future reference we fix some notations. Throughout this paper $k$ will denote a
totally real algebraic number field with the discriminant $D_k$, ring of integers $\cO$ and
ad\`ele ring $\A$. The set of places $V$ of $k$ is a union of the set $V_\infty$ of archimedean
and $V_f$ of finite places.
For $v\in V_f$, as usually, $k_v$ denotes the completion of $k$ at $v$, $\cO_v$ is the ring
of integers of $k_v$ with the uniformizer $\pi_v$ and the residue degree $\#\cO_v/\pi_v = q_v$.

\section{The volume formula}

\subsection{}
In a fundamental paper \cite{P} G.~Prasad obtained a formula for the volume of a principal
arithmetic subgroup of an arbitrary quasi-simple, simply connected group. This is an
extensive generalization of the results of Siegel, Tamagawa, Harder, Borel and other
people who worked in this direction. Gross has extended Prasad's formula to the
arithmetic subgroups of reductive groups~\cite{Gross}. From these results we can
write down a closed formula for the volume of a principal arithmetic subgroup $\Lambda$ of
$G = SO(1,2r)^o$.

We use the Euler-Poincar\'e normalization of the Haar measure on $G(\A)$
in the sense of Serre \cite{Serre}. Namely, for a discrete subgroup $\Lambda$
with finite covolume, we then have:
$$
|\chi(\Lambda\backslash G)| = \mu^{EP}(\Lambda\backslash G).
$$
For a principal arithmetic subgroup $\Lambda$ associated to a coherent collection of parahoric
subgroups $P = (P_v)_{v\in V_f}$:
$$
\mu^{EP}(\Lambda\backslash G) = \mu(\Lambda\backslash G) =
c_\infty D_k^{\frac{1}{2}\dim G}\left(\prod_{i=1}^{r}\frac{m_i!}{(2\pi)^{m_i+1}}\right)^{[k:\Q]}
\tau(G)\:\E\prod_{v\in T}\lambda_v,
\eqno(1)
$$
where
\begin{itemize}
\item[-] $c_\infty = 2$ is the Euler-Poincar\'e characteristic of the compact dual of
the symmetric space $G/SO(n)$ (see \cite{Serre}, \S3 and also \cite{BP}, \S4 for
the discussion);
\item[-] dimension $\dim G$ and exponents $\{m_i\}$ of our group $G$ of type $B_r$
are well known to be $\dim G = 2r^2 + r$ and $m_i = 2i-1\ (i = 1,\dots,r)$ \cite{Bourb};
\item[-] the Tamagawa number $\tau(G)=2$ (see \cite{Weil});
\item[-] $\E$ is an Euler product which in our case is given by
\ $\E = \zeta_k(2)\cdot\ldots\cdot\zeta_k(2r)$\ ($\zeta_k(.)$ is the Dedekind zeta
function of $k$);
\item[-] finally, the rational factors $\lambda_v\in\Q$ correspond to the (finite)
set $T$ of the finite places of $k$ over which $P_v \not\cong G_{qs}^o(\cO_v)$, where
$G_{qs}^o$ is the identity component of the quasi-split inner form of $G$.
\end{itemize}
This formula gives us the (generalized) Euler characteristic of $\Lambda$. The hyperbolic
volume of $\Lambda\backslash\Hy^n$ can be obtained from $|\chi(\Lambda\backslash G)|$ by
multiplying by the half of the volume of the unit sphere  $S^n$ in $R^{n+1}$:
$$\vol(\Lambda\backslash\Hy^{2r})  =
\frac{(2\pi)^r}{1\cdot 3 \cdot\ldots\cdot (2r-1)}\cdot |\chi(\Lambda\backslash G)|.$$
Note, that in odd dimensions the Euler characteristic vanishes  but we can still obtain
a similar formula for the covolume of an arithmetic subgroup without passing through the
Euler-Poincar\'e measure.

\subsection{} The $\lambda$-factors in $(1)$ are the most subtle matter. Fortunately, we can
explicitly compute the factors using the Bruhat-Tits theory. In \cite{GHY} this was done
for the parahoric subgroups which arise as the stabilizers of the maximal lattices. We will
extend the table from \cite{GHY} for the odd special orthogonal groups to the other maximal
parahorics.

Consider orthogonal group $G = SO_{2r+1}$ over a nonarchimedean local field $k_v$ whose residue
field $\cO/\pi$ has order $q$. By~\cite{Tits} the group $G$ belongs to one of the two possible
types: $B_r$ or $^2B_r$. For $r>2$ the local Dynkin diagrams and relative local index for the
nonsplit type are given on Figure~1. For $r=2$ there is an isogeny between the groups of types
$B$ and $C$, in~\cite{Tits} this case is represented by the diagrams of $C_2$ and $^2C_2$. We
leave to the reader to check that (with the suitable notations) all our computations remain to
be valid for this case as well.

\begin{figure}[ht]
\psfig{file=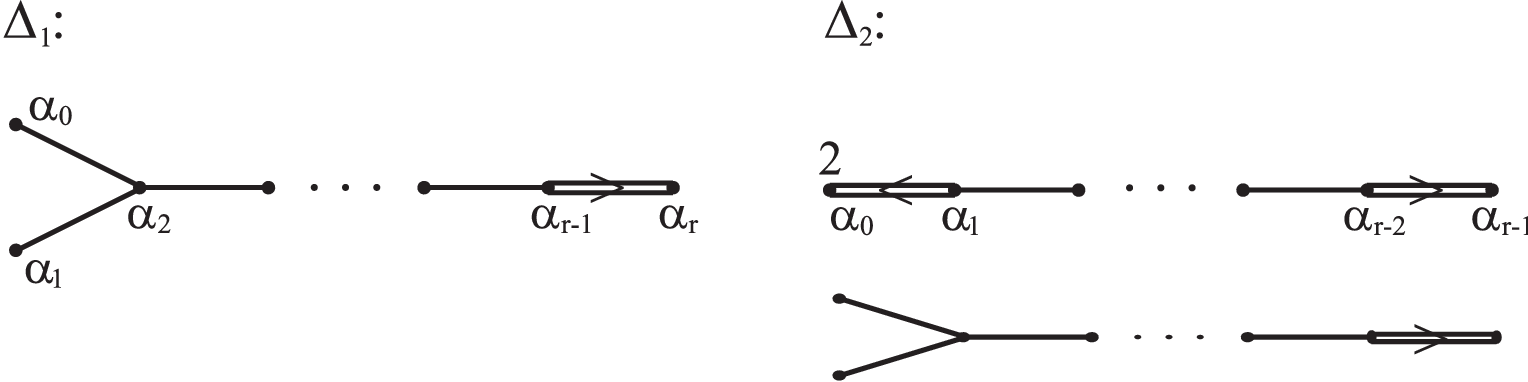, scale=.75}
\caption{}
\end{figure}

Similarly to~\cite{CR}, having the local diagrams we can enumerate all the types of the maximal
parahoric subgroups $P$ of $G$, and the type defines a parahoric subgroup up to conjugation
in $G$. However, some parahoric subgroups which are not conjugate in the
simply connected group can become conjugate in the adjoint group. This happens exactly when the
diagrams defining the types of the parahoric subgroups are symmetric with respect to an
automorphism of the entire diagram. So, in our case types $\Delta_1\backslash\{\alpha_0\}$
and $\Delta_1\backslash\{\alpha_1\}$ define conjugate subgroups.

For each type, using results of Bruhat and Tits (\cite{BT}, the
account of what we need can be found in~\cite{Tits}, \S3), we can determine the type
of the maximal reductive quotient ${\overline G}$ of the special fiber $\underline G$ of the
Bruhat-Tits group scheme associated with $P$ and also the type of the reductive quotient
${\overline G}_{qs}^o$ of the smooth affine group scheme ${\underline G}_{qs}^o$ which was defined
in ~\cite{Gross},~\S4 for the quasi-split inner form of $G$ (see also~\cite{P}).
Now, using the tables of orders of finite groups of Lie type (e.g. \cite{Ono}, Table~1)
for each of the cases the corresponding $\lambda$-factor is readily computed
by the formula from~\cite{GHY},~\S2:
$$
\lambda = \lambda (P)= \frac{q^{-N({\overline G}_{qs}^o)}\cdot \#{\overline G}_{qs}^o(\cO/\pi)}
{q^{-N({\overline G})}\cdot \#{\overline G}(\cO/\pi)}
$$
($N(\overline{G})$ denotes the number of positive roots of $\overline G$ over the
algebraic closure of the residue field $\cO/\pi$). We list the results in the following
table.

$$
\def\arraystretch{1.6}
\arraycolsep=4pt
\begin{array}{|c|c|c|c|}
\hline
{\rm maximal\ type}\ \theta_v & \overline G & {\overline G}_{qs}^o & \lambda\\
\hline
\Delta_1\backslash\{\alpha_0\}  & SO_{2r+1} & SO_{2r+1} & 1\\
\Delta_1\backslash\{\alpha_0, \alpha_1\} & GL_1\times SO_{2r-1} & SO_{2r+1} & \frac{q^{2r}-1}{q-1} \\
\Delta_1\backslash\{\alpha_i\},\ i=2,...,r-1 & O_{2i}\times SO_{2(r-i)+1} & SO_{2r+1} &
\frac{(q^i+1)\prod_{\nu=i+1}^r(q^{2\nu}-1)}{2\cdot\prod_{\nu=1}^{r-i}(q^{2\nu}-1)} \\
\Delta_1\backslash\{\alpha_r\} & O_{2r} & SO_{2r+1} & \frac{q^r+1}{2} \\
\Delta_2\backslash\{\alpha_0\} & ^2O_2\times SO_{2r-1} & SO_{2r+1} & \frac{q^{2r}-1}{2(q+1)} \\
\Delta_2\backslash\{\alpha_i\},\ i=1,...,r-2 & ^2O_{2(i+1)}\times SO_{2(r-i)-1} & SO_{2r+1} &
\frac{(q^{i+1}-1)\prod_{\nu=i+2}^r(q^{2\nu}-1)}{2\cdot\prod_{\nu=1}^{r-i-1}(q^{2\nu}-1)} \\
\Delta_2\backslash\{\alpha_{r-1}\} & ^2O_{2r} & SO_{2r+1} & \frac{q^r-1}{2} \\
\hline
\end{array}
$$

\medskip
\begin{center} {\sc Table 1} \end{center}
\medskip

\subsection{}
\begin{prop}
For any rank $r\ge 2$ and  $P_v \not\cong {\underline G}_{qs}^o(\cO_v)$ we have:
\begin{itemize}
\item[1)] $\lambda_v > 1$;
\item[2)] $\lambda_v > 2$ except for the case $r=2$, $q_v =2$, $\theta_v = \Delta_2\backslash\{\alpha_1\}$
\end{itemize}
\end{prop}
\begin{proof}
If $P_v$ is a maximal parahoric subgroup then the statement reduces to an easy
check of the values of $\lambda$ in Table~1. For an arbitrary parahoric subgroup
$P_v\subset G(k_v)$ there exists a maximal parahoric $P$ which contains $P_v$.
By the formula for $\lambda$ we get $\lambda(P_v) = [P:P_v]\:\lambda(P)$, so the
$\lambda$-factors of non-maximal parahoric subgroups also satisfy the conditions
(1) and (2).
\end{proof}

We remark that as in the simply connected case (\cite{BP}, Section~3.1 and Appendix~A)
the minimal values of $\lambda_v$ correspond to the special parahoric subgroups,
i.e. those for which the diagram representing $\theta_v$ is the Coxeter diagram of
the underlying finite reflection group.

\subsection{} \label{example}
{\it Example.} Let us consider the case when $f$ is the unimodular
integral quadratic form $-x_0^2 + x_1^2 + \dots + x_n^2$ ($n$ is even) and $\Gamma$
is the stabilizer of a maximal lattice $L$ on which $f$ takes integral values. Then
by~\cite{GHY} the group $\Gamma$ is a principal arithmetic subgroup of
$SO(f)^o = SO(1,n)^o$ and the types of the corresponding parahoric subgroups
can be determined from the local invariants of the quadratic form $f$. We have
$k = \Q$, $\cO = \Z$, $v$ runs through the primes of $\Q$, the determinant
$\delta_v(f) = \pm 1 \in \Q^\times$, the Hasse-Witt invariant $w_v(f) = \epsilon (SO(f;k_v)^o)$
is $1$ for $v\neq 2$ and $w_2(f) = 1$ if $n\equiv 0,2\ ({\rm mod}\ 8)$,
$w_2(f) = -1$ if $n\equiv 4,6\ ({\rm mod}\ 8)$. So we take the values for
$\lambda_v$ from the Table~4 in~\cite{GHY} and immediately obtain:
$$|\chi(\Gamma)| = 2\mkern-4mu\prod_{i=1}^{r}\frac{(2i-1)!}{(2\pi)^{2i}}
\ 2\mkern-4mu\prod_{i=1}^{r}\zeta(2i)\lambda_2 = 4\mkern-4mu\prod_{i=1}^r\frac{|B_{2i}|}{4i}
\cdot\left\{\begin{array}{@{\mkern-1mu}l@{\ {\rm if}\ }l}
1 & r\mkern-4mu\equiv 0,1 ({\rm mod}\ 4),\\
6^{-1}(2^{2r}-1) & r\mkern-4mu\equiv 2,3 ({\rm mod}\ 4)
\end{array}
\right.
$$
($B_{2i}$ are Bernoulli numbers: $B_2 = 1/6,\ B_4 = -1/30,\ B_6 = 1/42$ \dots).
This can be compared with \cite{RT} where the authors evaluated the
Siegel's limit, but the results will not coincide. The reason is that Ratcliffe and
Tschantz consider the arithmetic subgroups  $SO(1,n;\Z)$ of $G$ which are
the stabilizers of not maximal but unimodular lattices in the corresponding
quadratic spaces. The relation between these two cases is not straightforward,
but it appears that it is still possible to use a similar approach to obtain
the covolumes of the stabilizers of the unimodular lattices and, in particular,
to deduce the results of~\cite{RT}. We will explain this in detail in~\cite{BG}.

\subsection{}\label{max2}
We now return to the maximal arithmetic subgroups. As it was already mentioned
in Section~\ref{max1} any maximal arithmetic subgroup $\Gamma$ of $G$ can be
obtained as a normalizer in $G$ of a principal arithmetic subgroup. So, in order
to have a control over the volumes of $\Gamma = N_G(\Lambda)$ we need an estimate
for the index $[\Gamma : \Lambda]$. Following~\cite{BP} such an estimate can be
obtained from an exact sequence for the Galois cohomology of $k$ due to
Rohlfs~\cite{Rohlfs}. We have (\cite{BP}, Section~2.10):
$$
[\Gamma : \Lambda] \leq \#\prod_{v\in{\scS}}C(k_v)\cdot\#H^1(k,\widetilde C)_\xi
\cdot\mkern-8mu\prod_{v\in V\backslash S}\#\Xi_{\theta_v},
$$
where $C$ is the center of $G$, $\widetilde C$ is the center of its simply
connected inner form, $\Xi_{\theta_v}$ is the subgroup of the group of
automorphisms of the local Dynkin diagram of $G(k_v)$ stabilizing the type
$\theta_v$ and all the other notations can be found in~\cite{BP}.
In our case the center of $G$ is trivial, $S = V_\infty$, so we get
$$
[\Gamma : \Lambda] \le \#H^1(k,\widetilde C)_\xi
\cdot\mkern-8mu\prod_{v\in V_f}\#\Xi_{\theta_v}.
\eqno(2)
$$
Let $T_{ns}$ denote the (finite) set of places of $k$ for which $G$ does not split
over $k_v$. By~\cite{BP}, Proposition~5.1 applied to our group
$$
\#H^1(k,\widetilde C)_\xi \le h_k\cdot 2^{[k:\Q]+\#T_{ns}}.
\eqno(3)
$$
By~\cite{BP}, the proof of Proposition~6.1, the class number $h_k$ can be
estimated as
$$
h_k \le 10^2\left(\frac{\pi}{12}\right)^{[k:\Q]}D_k
\eqno(4)
$$
(this bound follows from the Brauer-Siegel theorem and Zimmert's bound for the
regulator of $k$. We refer to~\cite{BP} for the details.)
\medskip

We obtain: 
\begin{eqnarray*}
\lefteqn{\mu(\Gamma\backslash G)} &&\\
& \geq &
\frac{1}{h_k 2^{[k:\Q]+\#T_{ns}}\prod_{v\in V_f}\#\Xi_{\theta_v}}
\:c_\infty D_k^{\frac{1}{2}\dim G}\left(\prod_{i=1}^{r}\frac{m_i!}{(2\pi)^{m_i+1}}\right)^{[k:\Q]}
\mkern-16mu\tau(G)\:\E\prod_{v\in T}\lambda_v \\
& \geq & \frac{4 D_k^{\frac{1}{2}\dim G - 1}}{10^2 \left(\frac{\pi}{6}\right)^{[k:\Q]}}
\left(\prod_{i=1}^{r}\frac{m_i!}{(2\pi)^{m_i+1}}\right)^{[k:\Q]}\mkern-16mu\E
\prod_{v\in T}\lambda_v \left(\prod_{v\in V_f}\#\Xi_{\theta_v}\right)^{\mkern-12mu -1} \mkern-10mu 2^{-\#T_{ns}}.
\hspace*{3em} (5)
\end{eqnarray*}

\noindent
The group $\Xi_{\theta_v}$ is trivial if $G(k_v)$ is nonsplit or $\theta_v = \Delta_1\backslash\{\alpha_0\}$
and has order at most $2$ in all the rest of the cases, so we always have
$$
\prod_{v\in T}\lambda_v \left(\prod_{v\in V_f}\#\Xi_{\theta_v}\right)^{\mkern-12mu -1} \mkern-10mu 2^{-\#T_{ns}}
\geq \prod_{v\in T}\frac{\lambda_v}{2}.
$$

\subsection{}
\begin{prop}
$$
\E\prod_{v\in T}\lambda_v \left(\prod_{v\in V_f}\#\Xi_{\theta_v}\right)^{\mkern-12mu -1} \mkern-10mu 2^{-\#T_{ns}} > 1.
$$
\end{prop}
\begin{proof}
Except for the case $r=2$ and there exist $v\in V(k)$ with $q_v = 2$ the statement immediately follows from
Proposition~3.3~(2). For the remaining case we need to split the factors corresponding to the $2$-adic places
from the Euler product $\zeta_k(2)$. The meaning of this is that in order to make the estimate we need to
consider (almost) the actual volumes of certain parahoric subgroups, not just their quotients by the volume
of the standard parahoric which are captured in the $\lambda$-factors:
$$
\zeta_k(2)\zeta_k(4) \prod_{\substack{v\in T \\ q_v = 2 \\ }}\frac{\lambda_v}{2}
\prod_{\substack{v\in T \\ q_v\neq 2}}\frac{\lambda_v}{2} >
\zeta_k(2) \prod_{\substack{v\in T \\ q_v = 2}}\frac{\lambda_v}{2} \ge
\zeta'_k(2) \prod_{\substack{v\in T \\ q_v = 2}}\frac{1}{1-2^{-2}} \frac{2^2-1}{4} > 1.
$$
\end{proof}

\section{Orbifolds}

\subsection{} \label{thm1} \begin{theorem}
For any $n = 2r\ge 4$ there exists a {\it unique} compact orientable arithmetic hyperbolic $n$-orbifold
$O_{min}^n$ of the smallest volume. It is defined over the field $\Q[\sqrt{5}]$
and has Euler characteristic
$$
|\chi(O_{min}^n)| = \frac{\lambda(r)}{N(r)4^{r-1}} \prod_{i=1}^{r}|\zeta_{\Q[\sqrt{5}]}(1-2i)|,
$$
where: $\lambda(r) = 1$ if $r$ is even and $\lambda(r) =2^{-1}(4^r-1)$ if $r$ is odd;\\
\phantom{where:} $N(r)$ is a positive integer, $\le 4$ if $r$ is even, and $\le 8$
if $r$ is odd.
\end{theorem}

\begin{proof} {\bf 1.}
We are looking for a maximal arithmetic subgroup of $G = SO(1,n)^o$
of the smallest volume which is defined over $k\neq\Q$.

Let $k = \Q[\sqrt{5}]$. Then $k$ is a totally real quadratic field of the smallest
discriminant $D_k = 5$. By $(1)$ the volume is proportional to ${D_k}^{\dim G/2}$
and depends exponentially on the degree of the field (for big enough $r$),
so the smallest orbifold $O_0^n = \Gamma_0^n\backslash\Hy^n$ defined over $k$ is
a good candidate for $O_{min}^n$. Let $P = (P_v)_{v\in V_f(k)}$ be a coherent
collection of parahoric subgroups of $G(k)$ such that:
\begin{itemize}
\item[-] if $r$ is even $P_v = G(\cO_v)$ for all $v\in V_f$;
\item[-] if $r$ is odd $P_v = G(\cO_v)$ for all $v$ except one with the residue
characteristic $2$, for the remaining place $v_2$ we choose $P_{v_2}$ so that
$\lambda_{v_2} = (q^r-1)/2$.
\end{itemize}
These collections of parahoric subgroups satisfy the conditions of the maximality
criterion~\cite{CR} and the parity condition (Section~\ref{par_cond}).
Let
$$\Lambda_0^n = G(k)\cap\prod_v P_v.$$
So $\Lambda_0^n$ is a principal arithmetic subgroup of $G$. Let
$$\Gamma_0^n = N_G(\Lambda_0^n).$$
Then $\Gamma_0^n$ is a maximal arithmetic subgroup. We have:
\begin{eqnarray*}
\mu(\Lambda_0^n\backslash G) & = & 4 \cdot 5^{r^2+r/2}C(r)^2 \prod_{i=1}^{r}\zeta_k(2i)\;\lambda(r),\\
\mu(\Gamma_0^n\backslash G) & = & \frac{\mu(\Lambda_0^n\backslash G)}{[\Gamma_0^n:\Lambda_0^n]} =
\frac{\mu(\Lambda_0^n\backslash G)}{N(r)}.
\end{eqnarray*}
Here $C(r)$ denotes the product $\prod_{i=1}^{r}(2i-1)!/(2\pi)^{2i}$, $\lambda(r) = \lambda_{v_2}$
is as in the statement of the theorem and $N(r)$ is the order of the group of outer
automorphisms of $\Lambda_0^n$. By $(2)$ and $(3)$, $N(r)\le 4$ for even $r$ and $N(r)\le 8$
for odd $r$. Later on we will state a conjecture about the actual value of $N(r)$.

Consider the other groups defined over quadratic extensions of $\Q$. Note, that for odd
$r$ the set $T$ (of ``bad places'') should contain at least one place due to the parity
condition. For $\Lambda_0^n$ we have chosen $T$ in such a way that the $\lambda$-factor
in the volume formula has the smallest possible value for the groups defined over
$\Q[\sqrt{5}]$. So let $\Gamma^n$ be a maximal arithmetic subgroup of $G$ defined over
a totally real quadratic field $k$, $k \neq \Q[\sqrt{5}]$. By inequality $(5)$ and
Proposition~3.6:
$$ \mu(\Gamma^n\backslash G) > \frac{1}{h_k}\cdot D_k^{r^2+r/2} C(r)^2. $$
Now, except for the case $k = \Q[\sqrt{2}]$, $r = 3$:
$$
\frac{1}{h_k} D_k^{r^2+r/2} C(r)^2
\ge 4\cdot 5 ^{r^2+r/2} C(r)^2 2 \lambda(r) \ge
$$
$$
4\cdot 5 ^{r^2+r/2} C(r)^2 \prod_{i=1}^r\zeta_{\Q[\sqrt{5}]}(2i)\lambda(r)
\ge \mu(\Gamma_0^n\backslash G).
$$
In the first inequality for $D_k > 28$ we used the bound $(4)$ for the class number
$h_k$, for the remaining fields of the small discriminants the class numbers are
known to be equal to $1$. The second inequality is provided by the following
property of $\zeta_{\Q[\sqrt{5}]}(s)$:
$$
\prod_{i=1}^r\zeta_{\Q[\sqrt{5}]}(2i) < 2 \ \ {\rm for\ any}\ r. \eqno(*)
$$
\bigskip

\noindent
The proof of $(*)$ is easy. Let again $k = \Q[\sqrt{5}]$.
$$
P := \prod_{i=1}^r\zeta_k(2i) \leq \prod_{i=1}^\infty\zeta_k(2i)
= \zeta_k(2)\prod_{i=2}^\infty\zeta_k(2i) \leq \zeta_k(2)\prod_{i=2}^\infty\zeta^2(2i),
$$
$\zeta(n) = 1 + 1/2^n + 1/3^n + \ldots$\ is the Riemann zeta function.
By induction on $n$, for $n\ge 4$ we have $\zeta(n) \leq 1 + 2/2^n$.
So
$$
P \leq \zeta_k(2)\prod_{i=2}^\infty(1+ 2/2^{2i})^2.
$$
The right-side product converge and all the factors are $>1$ so we can take its logarithm:
$$
\log\left(\prod_{i=2}^\infty(1 + 2/2^{2i})^2\right)
=  \sum_{i=2}^\infty 2\log(1+1/2^{2i-1}) < \sum_{i=2}^\infty 1/2^{2i-2} = 1/3;
$$
$$
P < \zeta_k(2) e^{1/3} < 2
$$
and $(*)$ is proved.
\bigskip

The remaining case $k = \Q[\sqrt{2}]$, $r = 3$ is checked directly. We use Proposition~3.3
to estimate the $\lambda$-factors keeping the Euler product for the next inequality:
$$
\mu(\Gamma^n\backslash G) \ge
8^{10,5} C(r)^2 \prod_{i=1}^3\zeta_k(2i)\prod_{v\in T} \frac{\lambda_v}{2} \ge
$$
$$
8^{10,5} C(r)^2 \prod_{i=1}^3\zeta_k(2i) >
4\cdot5^{10.5} C(r)^2 \prod_{i=1}^3\zeta_{\Q[\sqrt{5}]}(2i)\frac{4^3-1}{2}.
$$
So we are done with the quadratic fields and can proceed to the higher degrees.
\medskip

Let $[k:\Q] = 3$. We have:
\begin{eqnarray*}
\frac{\mu(\Gamma^n\backslash G)}{\mu(\Gamma_0^n\backslash G)}
& \geq &
\left(\frac{D_k}{5}\right)^{r^2+r/2} \frac{C(r)\prod_{i=1}^{r}\zeta_k(2i)\prod_{v\in T}\frac{\lambda_v}{2}}
{2^3h_k \prod_{i=1}^{r}\zeta_{\Q[\sqrt{5}]}(2i)\;\lambda(r)}.
\end{eqnarray*}


First consider the totally real cubic field of the smallest discriminant $D_k=49$. This field
has $h_k = 1$, moreover, since its ring of integers does not have prime ideals of norm $2$, we
can use Proposition~3.3 to estimate the $\lambda$-factors in all the cases:
$$
\frac{\mu(\Gamma^n\backslash G)}{\mu(\Gamma_0^n\backslash G)} >
\left(\frac{49}{5}\right)^{r^2+r/2}
\frac{C(r)\prod_{i=1}^{r}\zeta_k(2i)}{8\prod_{i=1}^{r}\zeta_{\Q[\sqrt{5}]}(2i)\;\lambda(r)}
> 1.
$$
This inequality can be checked directly for $r =2$; for the higher ranks it is enough
to estimate the product of $\zeta_k(2i)$ by $1$ from below and the product of
$\zeta_{\Q[\sqrt{5}]}(2i)$ by~$(*)$ from above, which gives an easy-to-check inequality.

For the other cubic fields by Proposition~3.6 and inequality $(*)$ we have
$$
\frac{\mu(\Gamma^n\backslash G)}{\mu(\Gamma_0^n\backslash G)} >
\left(\frac{D_k}{5}\right)^{r^2+r/2} \frac{C(r)}
{8 h_k 2 \lambda(r)}.
$$
Again, using the precise values of $h_k$ for the fields of the small discriminants ($D_k = 81$, $148$,
$169$) and bound~(4) for the other fields, we see that this is always greater then $1$.
\medskip

For $d = [k:\Q] \geq 4$ we will make use of the known lower bounds for the discriminants of
the totally real number fields (see \cite{Odlyzko}):
\begin{itemize}
\item[] if $d = 4$\ \ $D_k > 5^d$;
\item[] if $d = 5$\ \ $D_k > 6.5^d$;
\item[] if $d\geq 6$\ \ $D_k > 7.9^d$ (and $D_k > 10^d$ if $d\geq 8$).
\end{itemize}
The cases $d = 4,\ 5$ are considered similar to the previous case $d = 3$ and we allow ourselves
to skip the details. Let $\Gamma^n$ is defined over a field $k$ of degree $d\geq 6$. We have:
\begin{eqnarray*}
\frac{\mu(\Gamma^n\backslash G)}{\mu(\Gamma_0^n\backslash G)}
& > & \frac15\left(\frac{D_k}{5}\right)^{r^2+r/2-1}
\frac{C(r)^{d-2}}{2^d\cdot 100\cdot\left(\frac{\pi}{12}\right)^d 2\;\lambda(r)}.
\end{eqnarray*}
For $r\geq 3$ we can estimate $D_k$ by $7.9^d$ and then show that
$\mu(\Gamma^n\backslash G)/\mu(\Gamma_0^n\backslash G) > 1$ for any $d$. This does not work for
$r = 2$ since for large $d$ the factor $C(r)^{d-2}$ becomes too small. In order to get rid of it
we use the second bound $D_k > 10^d$ for $d\geq 8$. Note, that here we do not need any
particular knowledge of the class numbers.
\medskip

We proved that $O_{min}^n = O_0^n = \Gamma_0^n\backslash\Hy^n$ has the smallest possible volume
for each $n$. Using the functional equation for the Dedekind zeta function we can write down
the formula for $\mu(\Gamma_0^n\backslash G) = |\chi(O_{min}^n)|$ in a compact form which is given
in the statement. It remains to show the uniqueness of $O_{min}^n$.
\medskip

\noindent{\bf 2.}
Let $H_1/k$ and $H_2/k$ be two algebraic groups defined over $k = \Q[\sqrt5]$ such that each
$H_i(k\otimes\R)$ admits a surjective homomorphism onto $G$ with a compact kernel. Then $H_1$ and
$H_2$ are $k$-isogenous and the isogeny takes arithmetic subgroups to arithmetic subgroups. So we
can fix an algebraic group $H/k$ and the surjective homomorphism with a compact kernel
$\phi: H(k\otimes\R)\to G$, such that $H$ is of type $B_r$ and can be supposed to be the adjoint group
since $G$ is centerless. Let $\Lambda_1$ and $\Lambda_2$ be two arithmetic subgroups of $H(k)$
of the same maximal type $P = (P_v)$. We want to prove that $N_G(\phi(\Lambda_1))$ is conjugate in $G$
to $N_G(\phi(\Lambda_2))$.

For each finite place $v\in V_f$ there exists $g_v$ such that $\Lambda_{1,v} =
g_v \Lambda_{2,v} g_v^{-1}$ (see Section~3.2). The set of places
$$ S = \{ v\in V_f \mid \Lambda_{1,v}\neq H(\cO_v) {\rm\ and\ } \Lambda_{2,v}\neq H(\cO_v)\} $$
is finite. This is a known fact. To prove it one can first show the finiteness of such a
set of places for the simply connected inner form of $H$ using the strong approximation
property (see e.g. \cite{CR}), and then transfer it to $H$ itself by an inner twist. So the set
$$ U = \{ g\in H(\A_f) \mid \prod_{v\in V_f} \Lambda_{1,v} = g (\prod_{v\in V_f}\Lambda_{2,v})
g^{-1} \} $$
is not empty. Moreover, for each $g = (g_v) \in U$ we have an open subset
$$\prod_{v\in S} g_v\Lambda_{2,v} \prod_{v\in V_f\backslash S} \Lambda_{2,v} \subset U.$$

If $S = \emptyset$ then $\Lambda_1 = \Lambda_2$ and there is nothing to prove. Suppose
$S$ is non empty. We consider $H(k)$ diagonally embedded into $H(\A_f)$ and in $\prod_{v\in S} H(k_v)$.
By the weak approximation property for $H$ ($H$ is an adjoint group
so it has weak approximation \cite{Harder}), $H(k)\cap\prod_{v\in S} g_v\Lambda_{2,v}$
is dense in $\prod_{v\in S} g_v\Lambda_{2,v}$ with respect to the product topology,
so there exists a non empty open subset $X\subset H(k)$ such that for any $x\in X$ and any
$v\in S$:
$$\Lambda_{1,v} = x\Lambda_{2,v}x^{-1}.$$
Let $p_1,\dots ,\ p_n$ be the set of prime ideals in $\cO$
which define the places from $S$. Consider the ring $R = \cO[\frac{1}{p_1}, \dots ,\ \frac{1}{p_n}]$.
Since $S\neq\emptyset$ it is a dense subset of $k$, and so $H(R)$ is dense in $H(k)$.
Consequently there exists $r\in X\cap H(R)$. We have:
$$
r \Lambda_2 r^{-1} = r ( H(k)\cap\prod_{v\in V_f}\Lambda_{2,v} ) r^{-1} =
H(k)\cap r ( \prod_{v\in V_f}\Lambda_{2,v} ) r^{-1} =
$$
$$
H(k)\cap \prod_{v\in S}r\Lambda_{2,v}r^{-1} \prod_{v\in V_f\backslash S}r H(\cO_v) r^{-1} =
H(k)\cap \prod_{v\in V_f}\Lambda_{1,v} = \Lambda_1;
$$

$$ \phi(r) \phi(\Lambda_2) \phi(r)^{-1} = \phi(r\Lambda_2r^{-1}) = \phi(\Lambda_1); $$

$$ \phi(r) N_G(\phi(\Lambda_2)) \phi(r)^{-1} = N_G(\phi(\Lambda_1)).$$

\end{proof}
\medskip

The reader can notice that part~1 of the proof can be simplified if we suppose that
$r>2$ but the case $r=2$ is important. We will come back to it in the next section.
Now let us give some remarks concerning the general case.
\medskip

\subsection{} \label{conj}
The value of $N(r)$ is the order of the outer automorphisms group of the
principal arithmetic subgroup $\Lambda_0^n$. Since we are in a very extremal situation
we suppose that, in fact, $\Lambda_0^n$ has no non-trivial symmetries. This can be
checked for small $n$ for which the group $\Lambda_0^n$ is reflective (see \cite{V}).
We do not know how to prove this observation for the higher dimensions, but
still we would like to have it as a conjecture.
\medskip

\begin{conjecture}
For all $r\ge 2$,\ $N(r) = 1$ and so
$$
|\chi(O_{min}^n)| = \frac{\lambda(r)}{4^{r-1}}\prod_{i=1}^{r}|\zeta_{\Q[\sqrt5]}(1 - 2i)|.
$$
\end{conjecture}

\subsection{} \label{rem2} We can describe groups $\Gamma_0^n$ of the smallest orbifolds
as the stabilizers of lattices in quadratic spaces. Let
$$f = -\frac{1+\sqrt{5}}{2}x_0^2 + x_1^2 + \ldots + x_n^2.$$
and $(V,f)$ is the corresponding $(n+1)$-dimensional quadratic space. By~\cite{GHY}
for even $r = n/2$ the coherent collection of parahoric subgroups defining
the principal arithmetic subgroups $\Lambda_0^n$ has the same type as
the one that gives the stabilizer of the maximal lattice in $(V,f)$.
So by the uniqueness argument from the proof of the theorem,
$\Lambda_0^n$ is the stabilizer of the maximal lattice in $(V,f)$.
Similarly, by \cite{BG} for $r$ odd $\Lambda_0^n$ is the stabilizer of the odd unimodular
lattice in $(V,f)$ or, equivalently, it is the stabilizer of the maximal lattice in
$(V, 2f)$. Consequently, the groups $\Gamma_0^n$ are obtained as the normalizers in $G$ of
the stabilizers of the lattices.

\subsection{} \label{rem3}
For completeness, let us also consider the non-compact case which is easy because
the only possible field of definition is $k = \Q$. Similarly to the previous constructions
(see also Example~\ref{example}) we obtain:

\begin{quote}
For any $n=2r\ge 4$ there is an unique non-compact orientable arithmetic hyperbolic
$n$-orbifold ${O'}^n_{min}$ of the smallest volume. It has Euler characteristic
$$
|\chi({O'}_{min}^n)| = \frac{\lambda'(r)}{N'(r)2^{r-2}}\prod_{i=1}^{r}|\zeta(1-2i)|.
$$
where: $\lambda'(r) = 1$ if $r\equiv 0,1\ ({\rm mod}\ 4)$,\newline
\phantom{where:} $\lambda'(r) = 2^{-1}(2^r-1)$ if $r\equiv 2,3\ ({\rm mod}\ 4)$;\newline
\phantom{where:} $N'(r)$ is a positive integer, $\le2$ if $r\equiv 0,1\ ({\rm mod}\ 4)$,\newline
\phantom{where:} and $\le 4$ if $r\equiv 2,3\ ({\rm mod}\ 4)$.

\noindent
For $r\equiv 0,1\ ({\rm mod}\ 4)$ the group of ${O'}_{min}^n$ is the normalizer of
the stabilizer of the maximal lattice in quadratic space $(V,f)$ defined by
$$f = -x_0^2 + x_1^2 + \ldots + x_n^2;$$
and for $r\equiv 2,3\ ({\rm mod}\ 4)$ it is the normalizer of the stabilizer
of the odd unimodular lattice in $(V,f)$.
\end{quote}

\noindent
Conjecture~\ref{conj} also applies to ${O'}_{min}^n$ and says that $N'(r) = 1$ for
all $r$.

\subsection{}\label{rem4}
We will now compute the Euler characteristics of the smallest orbifolds for small $n$.
We will give the values for the principal arithmetic subgroups and then either
Conjecture~\ref{conj} is true or one should divide by the actual value of $N(r)$ in
order to obtain the Euler characteristic of the smallest orbifolds. In any case, since
$N(r)$ is bounded and always smaller then $8$, this will not change the qualitative
picture. We have:
$$
\def\arraystretch{1.5}
\begin{array}{|c|*{6}{|c}}
\hline
  n = 2r \geq 4            & 4 & 6 & 8 & 10 & 12 & \\
\hline
|\chi(\Lambda_{min}^n)|    & \frac{1}{7200} & \frac{67}{576000} & \frac{24187}{8709120000} &
\frac{309479461547}{3483648000000} & \frac{7939510008126649607}{3766102179840000000} & \dots \\
\hline
|\chi({\Lambda'}_{min}^n)| & \frac{1}{960} & \frac{1}{207360} & \frac{1}{348364800} &
\frac{1}{91968307200} & \frac{691}{191294078976000} & \\
\hline
\end{array}
$$

$$
\def\arraystretch{1.5}
\begin{array}{c*{3}{|c}|}
\hline
    & 14 & 16 & 18 \\
\hline
\dots & 8.1824...\cdot10^{10} & 3.3481...\cdot10^{16} & 1.7455...\cdot10^{34}\\
\hline
    & \frac{87757}{289236647411712000} & \frac{2499347}{2360171042879569920000} & \frac{109638854849}{67802993719844284661760000} \\
\hline
\end{array}
$$

\medskip
\begin{center} {\sc Table 2} \end{center}
\medskip

\noindent
The smallest non-compact orbifold for all $n$ (which is also the smallest among all the
arithmetic orbifolds) has dimension $n = 16$ and $\chi = \frac{2499347}{2360171042879569920000} = 1.0589...\cdot10^{-15}$.
The smallest compact orbifold has dimension $n = 8$ and $\chi = \frac{24187}{8709120000} = 0.00000277...$\ .
For $n = 4$ the volume of the smallest compact orbifold is less then that of the non-compact
one and for all bigger $n$ the non-compact orbifolds are smaller. After $n = 10$  in the
compact case and $n = 18$ for non-compact (and so for all the arithmetic hyperbolic orbifolds)
the minimal volumes start to increase and then grow exponentially with respect to the
dimension.

It was first discovered in \cite{Lub} that the minimal covolume can be attained on a non-uniform
(that is, not cocompact) lattice. The result was obtained for the groups $SL(2,K)$ over local
fields $K$ of a positive characteristic. The natural question which appeared the same time was
whether this is a purely local phenomenon or it is also possible for the groups over global fields.
Our computation gives the answer to this question for the odd orthogonal groups over
the totally real number fields, moreover, the method indicates that the minimality of the covolume
of the non-uniform lattices might be always the case for the groups of a high enough rank.

\subsection{}\label{rem5}
The previous remark can be considered in a wider context of \cite{BP} where the discreteness
of the set of covolumes of arithmetic subgroups was proved in a very general setting.
In particular, it follows from \cite{BP} that there exist ``absolutely smallest'' among all
the $S$-arithmetic subgroups of $G$ over $k$ when $G$ runs through the algebraic $k$-groups
of absolute rank $\ge 2$, the global field $k$ can be either a number or a function field,
$S$ is any finite subset of places of $k$ containing all the archimedean places and the Haar
measures are chosen in a consistent way. Our results imply that for the adjoint groups $G$ of type
$B_r$ and real rank $1$ over the totally real number fields the smallest arithmetic subgroup
is the unique smallest arithmetic non-cocompact subgroup of the group $G$ of rank $r = 8$.
It is not hard to generalize this results (except the uniqueness) to the $S$-arithmetic subgroups
and to the other forms of type $B_r$ over totally real fields. Our previous remark allows to
conjecture that in general the smallest group might be non-cocompact which significantly
reduces the number of possible candidates. Still the detailed study of this question lies
beyond the scope of this paper. We have to point out that the geometric or any other meaning of
the absolutely smallest group is completely mysterious, the only thing we know is that such a
group or groups exist.

\section{$4$-Manifolds}

\subsection{} Let us now consider the problem of determining the smallest arithmetic
manifolds. This is a much more difficult task. From the arithmetic point of view the
first difficulty is that we can not just estimate the Euler products but rather
we have to deal with their rational structure. So, obtaining the precise values of
the Euler characteristic for the groups of interest will be the first step. After
this one needs to study the low index subgroups lattice of the distinguished
groups and find the torsion-free subgroups. We will restrict our attention to the
hyperbolic dimension~$4$.

\subsection{} The smallest known example of a compact orientable hyperbolic $4$-manifold
was constructed by Davis~\cite{Davis}. It can be shown that the Davis manifold $M_D$ is
arithmetic and defined over the field $\Q[\sqrt 5]$~\cite{EM}. The Euler characteristic
$\chi(M_D)$ is equal to $26$. We will be looking for smaller examples, so we are
interested in the manifolds with $\chi < 26$. It is well-known that the Euler characteristic
of a compact orientable $4$-manifold is an even positive integer which gives us one more
natural restriction on $\chi$.

\subsection{}\label{lemma}
We start with refining some estimates from the proof of Theorem~\ref{thm1} for the case $r=2$ and
a larger bound for $\chi$.
\medskip

\begin{lemma} If an orientable compact arithmetic hyperbolic $4$-manifold defined over a
field $k$ has $\chi\le 24$ then one of the following possibilities hold:
\begin{itemize}
\item[(1)] $d = 2$, $D_k \le 362$;
\item[(2)] $d = 3$, $D_k \le 3104$;
\item[(3)] $d = 4$, $D_k \le 26574$;
\item[(4)] $d = 5$, $D_k \le 227481$;
\item[(5)] $d = 6$, $D_k \le 1947276$
\end{itemize}
where $d = [k:\Q]$ is the degree of $k$ and $D_k$ is its discriminant.
\end{lemma}
\begin{proof}
The group of an arithmetic manifold is a (torsion-free) subgroup of a maximal arithmetic
subgroup $\Gamma$, so we have $\chi(\Gamma) \leq 24$. From the other side
\begin{eqnarray*}
\chi(\Gamma)
& > & \frac{1}{2^d h_k}\cdot4 D_k^5 \left(\frac{6}{2^6\pi^6}\right)^d \\
& \geq & \frac{1}{2^d\cdot 10^2\cdot(\frac{\pi}{12})^d D_k}\cdot4 D_k^5 \left(\frac{6}{2^6\pi^6}\right)^d\\
& \geq & D_k^4\cdot\frac{1}{25}\cdot\left(\frac{6^2}{2^6\pi^7}\right)^d.\\
\end{eqnarray*}
And thus we get
$$D_k < \left( 24\cdot25\cdot\left(\frac{6^2}{2^6\pi^7}\right)^d\right)^{1/4}.$$
For $d\geq 8$ this upper bound becomes smaller then the lower bound $10^d$ for the
discriminant of a totally real field of degree $d$ from~\cite{Odlyzko}. For $d=7$
the precise smallest value of $D_k$ is known to be $(11.051...)^7$~\cite{Odlyzko}, and
it again appears to be bigger then our upper bound for $D_k$. So we are left with the
$5$ remaining values of $d$ and for each of them we compute the corresponding upper bound
for $D_k$ from the above inequality.
\end{proof}

\subsection{} Using the tables of number fields of low degree~\cite{BFPOD} we can perform
a more careful analysis of the groups over the fields which satisfy the conditions of
Lemma~\ref{lemma}. There are many fields which fit the conditions and we used a
simple program for {\em GP/PARI}\ \ to perform the calculations. We obtain that among
$109 + 98 + 182 + 45 + 32 = 466$ totally real fields which have discriminant in
one of the $5$ ranges only $21 + 12 + 12 + 2 = 47$ can actually admit the groups with
$\chi\leq 24$ if we use the precise values of the class numbers in the volume estimate.
For the remaining fields we compute the numerator $\nu$ of the Euler characteristic
$\chi$ of the smallest arithmetic group $\Gamma$ defined over the field
($\Gamma$'s correspond to the principal arithmetic subgroups for which all the
$\lambda$-factors in the volume formula are equal to $1$). Since when passing to
a finite index subgroup $\Gamma'$ of $\Gamma$ the number $\nu$ still divides the
numerator of the Euler characteristic $\chi(\Gamma')$, we can discard all the
groups with $\nu > 24$ or $\nu$ is odd and $> 12$. There are also several maximal
groups with non-trivial $\lambda$-factors that fit into our range and these need
to be checked in an entirely similar way. Finally, we are left with only two groups
$\Gamma_1$ and $\Gamma_2$ which are defined over $\Q[\sqrt5]$, $\Q[\sqrt2]$ and
have $\chi = 1/7200,\ 11/5760$, respectively. Group $\Gamma_1$ is the
group of the smallest arithmetic $4$-orbifold (see Section~4), $\Gamma_2$ is
the smallest group defined over $\Q[\sqrt2]$ (the same argument as in the proof of
Theorem~\ref{thm1} can be used to show that $\Gamma_2$ is defined uniquely up to
conjugations in $G$).
\medskip

Let us summarize the results.

\subsection{}
\begin{theorem}\label{thm2}
If there exists a compact orientable arithmetic hyperbolic \linebreak
$4$-manifold $M$ having $\chi(M) \leq 24$ then it satisfies one of the following conditions:

\begin{itemize}
\item[1)]
$M$ is defined over $\Q[\sqrt{5}]$ and has the form  $\Gamma_M\backslash\Hy^4$ with $\Gamma_M$
is a torsion-free subgroup of index $7200\chi(M)$ of the group $\Gamma_1$ of the smallest
arithmetic $4$-orbifold;
\item[2)]
$M$ has Euler characteristic $22$, is defined over $\Q[\sqrt{2}]$, and its group is a torsion-free
subgroup of index $11520$ of $\Gamma_2$ which is the smallest principal arithmetic subgroup
of $SO(1,4)^o$ defined over $\Q[\sqrt2]$.
\end{itemize}
\end{theorem}

\subsection{} This result reduces the problem of finding the smallest compact arithmetic
$4$-manifold to a computational problem: we need to search for the ``low'' index torsion-free
subgroups of the groups $\Gamma_1$ and $\Gamma_2$ defined above. The first step to
implement this in practice is to find good presentations for the maximal groups. For
the group $\Gamma_1$ this can be done by identifying it with the orientation-preserving
subgroup of a Coxeter group $\Gamma_1'$ which has the Coxeter diagram given on
Figure~2.

\begin{figure}[ht]
\psfig{file=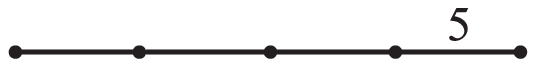, scale=.75}
\caption{}
\end{figure}
\noindent
(It is easy to check that $\Gamma_1'$ is an arithmetic subgroup of $O(1,4)$ defined
over $\Q[\sqrt5]$ and $\chi(\Gamma_1') = 1/14400$, so its orientation-preserving subgroup
is $\Gamma_1$ by the uniqueness of the smallest arithmetic orbifold.)

Now we can search for the torsion-free subgroups of $\Gamma_1'$ of index $14400\chi(M)$.
This, in principle, can be done by using the computer programs like {\em GAP}. The indexes
of the subgroups we are interested in are quite large, but, as it can be checked, the
Coxeter group $\Gamma_1'$ has not many subgroups of low index, and so the computation
looks more or less realistic.

We do not know whether or not the group $\Gamma_2$ is also reflective in a sense of~\cite{V},
and we suppose that it is not. Since we are dealing with a stabilizer of not a modular
but a maximal lattice the application of Vinberg's algorithm~\cite{V} for determining the
maximal subgroup generated by reflections is not straightforward here.
It is possible to  write down an explicit matrix representation for the generators of
$\Gamma_2$ in $SO(1,4)$, but we will not do it now. In any case this group can provide only
an example with $\chi = 22$ which is already quite large.

\end{document}